\documentstyle[10pt,amstex,amscd,graphicx]{amsart}

\numberwithin{equation}{section}

\newtheorem{thm}{Theorem}[section]

\newtheorem{prop}[thm]{Proposition}
\newtheorem{defn}[thm]{Definition}

\newcommand\AAA{{\mathcal A}}
\newcommand\BB{{\mathcal B}}

\newcommand\DD{{\mathcal D}}
\newcommand\EE{{\mathcal E}}
\newcommand\FF{{\mathcal F}}
\newcommand\GG{{\mathcal G}}
\newcommand\HH{{\mathcal H}}

\newcommand\LL{{\mathcal L}}
\newcommand\MM{{\mathcal M}}

\newcommand\PP{{\mathcal P}}

\newcommand\SSS{{\mathcal S}}
\newcommand\TT{{\mathcal T}}

\newcommand\PMF{{\PP\kern-2pt\MM\FF}}

\newcommand\PML{{\PP\kern-2pt\MM\LL}}

\newcommand\half{{\textstyle{\frac12}}}

\newcommand\hhat{\widehat}

\newcommand{\fsubd}{\mathrel{{\scriptstyle\searrow}\kern-1ex^d\kern0.5ex}}
\newcommand{\bsubd}{\mathrel{{\scriptstyle\swarrow}\kern-1.6ex^d\kern0.8ex}}
\newcommand{\fsubeq}{\mathrel{\raise-.7ex\hbox{$\overset{\searrow}{=}$}}}
\newcommand{\bsubeq}{\mathrel{\raise-.7ex\hbox{$\overset{\swarrow}{=}$}}}

\newcommand{\tsh}[1]{\left\{\kern-.9ex\left\{#1\right\}\kern-.9ex\right\}}

\begin{document}
\title{Notes on Stable Teichm\"{u}ller quasigeodesics}
\author{Abhijit Pal}
\address{Department of Mathematics and Statistics, IISER-Kolkata}

\date{}
\begin{abstract}
 In this note, we prove that for a cobounded,
Lipschitz path $\gamma:I\to\TT$
 if the pull back bundle $\mathcal H_{\gamma}$ over $I$ is a strongly relatively hyperbolic metric space
then there exists a geodesic $\xi$ in $\TT$ such that  $\gamma(I)$ and $\xi$ are close to each other.
\end{abstract}

\maketitle
Suppose $S_{g,n}$ is a  surface of genus $g$ with $n$ punctures such that its
Euler characteristic $\chi(S_{g,n})<0$.
Consider the Teichmuller space $\TT=Teich(S_{g,n})$ of $S_{g,n}$,  there is a  smooth fiber bundle $\SSS\to\TT$ over $\TT$,
whose fiber $\SSS_\sigma$ over $\sigma\in\TT$ is $S_{g,n}$ with metric $\sigma$. Let $\HH$ be the universal cover of $\SSS$,
then the universal covering $\HH\to\SSS$ defines a smooth fiber bundle $\HH\to\TT$ whose fiber $\HH_\sigma$ over $\sigma\in\TT$ is
isometric to the hyperbolic plane ${\mathbb{H}}^2$.
The purpose of this note is to prove that for a $\BB$-cobounded,
Lipschitz path $\gamma:I\to\TT$, where $\BB$ is a compact subset of $\TT$,
 if the pull back bundle $\mathcal H_{\gamma}$ over $I$ is a strongly relatively hyperbolic metric space
then there exists a geodesic $\xi$ in $\TT$ such that the Hausdorff distance between $\gamma(I)$ and $\xi$ 
is bounded.
This is a straightforward generalization of a result due to Mosher, Theorem 1.1 of \cite{mos}, where the statement was proven for
 closed surfaces admitting hyperbolic metrics with the assumption that $\HH_{\gamma}$ is a hyperbolic metric space.
\section{Relative Hyperbolicity}
Let $X$ be a path metric space. A collection of closed
 subsets $\DD = \{ D_\alpha\}$ of $X$ will be said to be {\bf uniformly
 separated} if there exists $\epsilon > 0$ such that
$d(D_1, D_2) \geq \epsilon$ for all distinct $D_1, D_2 \in \DD$.

\begin{defn} (Farb \cite{farb})
The {\bf electric space} (or coned-off space) $\EE(X,\DD)$
corresponding to the
pair $(X,\DD )$ is a metric space which consists of $X$ and a
collection of vertices $v_\alpha$ (one for each $D_\alpha \in \DD$)
such that each point of $D_\alpha$ is joined to (coned off at)
$v_\alpha$ by an edge of length $\half$. 
$X$ is said to be {\bf weakly hyperbolic} relative to the collection $\DD$ if $\EE(X,\DD)$ is a hyperbolic metric space.
\label{el-space}
\end{defn}
For a path $\gamma \subset X$, there is an induced path $\hhat{\gamma}$ in $\EE(X,\DD)$ obtained by coning the portions of
$\gamma$ lying in  sets $D\in\DD$. If $\hhat{\gamma}$ is a geodesic (resp. $P$-quasigeodesic) in $\EE(X,\DD)$, 
$\gamma$ is called a {\em relative geodesic} (resp.
{\em relative $P$-quasigeodesic}).
 \begin{defn}\cite{farb-relhyp}
Relative geodesics (resp. $P$-quasigeodesics) in  
$(X,\DD )$ are said to satisfy {\bf bounded region penetration properties} if there exists $K = K(P )>0$ such that
for any two  relative geodesics (resp. $P$-quasigeodesics without backtracking)
$\beta$, $\gamma$   joining $x, y \in X$ following two properties are satisfied:\\
(1)  if  precisely one of $\{ \beta , \gamma \}$ meets 
 a   set $D_\alpha$, 
then the length (measured in the intrinsic path-metric
  on  $D_\alpha$ ) from the first (entry) point
  to the last 
  (exit) point (of the relevant path) is at most $K$, \\
(2)  if both $\{ \beta , \gamma \}$ meet some  $D_\alpha $
 then the length (measured in the intrinsic path-metric
  on  $D_\alpha$ ) from the entry point of
 $\beta$ to that of $\gamma$ is at most $K$; similarly for exit points. \\
\end{defn}

\begin{defn} (Farb \cite{farb-relhyp} ) $X$ is said to be hyperbolic relative to the uniformly separated collection $\DD$ if 
 $X$ is weakly hyperbolic relative to $\DD$ and
 relative $P$ quasigeodesics without backtracking satisfy the bounded region penetration properties. \\
\end{defn}
\noindent {\underline{Gromov's definition of  relative hyperbolicity }:}\\

 \begin{defn}\cite{pal}
For any geodesic metric space 
$(D,d)$, the {\em hyperbolic cone} (analog of a horoball)
$D^h$ is the metric space 
$D\times [0,\infty) = D^h$ equipped with the 
path metric $d_h$ obtained from two pieces of
 data \\ 
1) $d_{h,t}((x,t),(y,t)) = 2^{-t}d_D(x,y)$, where $d_{h,t}$ is the induced path
metric on $D_t=D\times \{t\}$.  Paths joining 
$(x,t),(y,t)$ and lying on  $D_t=D\times \{t\}$ 
are called {\em horizontal paths}. \\
2) $d_h((x,t),(x,s))=\vert t-s \vert$ for all $x\in D$ and for all $t,s\in [0,\infty)$, and the corresponding paths are called 
{\em vertical paths}. \\
3)  for all $x,y \in D^h$,  $d_h(x,y)$ is the path metric induced by the collection of horizontal and vertical paths. \\
\end{defn}

\begin{defn}\label{gro}\cite{pal}
Let $\delta\geq 0$. 
Let $X$ be a geodesic metric space and $\DD$ be a collection of mutually disjoint uniformly separated subsets of $X$.  
$X$ is said to be  $\delta$-hyperbolic relative to $\DD$ in the sense of Gromov, 
if the quotient space $\GG(X,\DD)$,  obtained by attaching the hyperbolic cones
$ D^h$ to $D \in \DD$  via the identification $(x,0)\sim x$
  for all $x\in D$,  is a  $\delta$-hyperbolic metric space.
$X$ is said to be  hyperbolic relative to $\DD$ in the sense of Gromov if $\GG(X,\DD)$ is a $\delta$-hyperbolic
metric space for some $\delta\geq 0$.
\end{defn}
\begin{thm} (Bowditch \cite{bow})
Let $X$ be a geodesic metric space and $\DD$ be a collection of mutually disjoint uniformly separated subsets of $X$. 
 $X$ is  hyperbolic relative to the
collection $\DD$ of uniformly separated subsets of $X$ in the sense of Farb if and only if  $X$ is  hyperbolic relative to the
collection $\DD$ of uniformly separated subsets of $X$ in the sense of Gromov.
\end{thm}
\section{Main Theorem}
Suppose $p_1,...,p_n$ are the punctures of $S_{g,n}$, then  each Teichmuller metric $\sigma$ on $S_{g,n}$ corresponds to collections 
$\DD_{\sigma}(p_1),...,\DD_{\sigma}(p_n)$ of horodisks in the fiber $\HH_\sigma$ of the bundle 
$\HH\to\TT$ satisfying the following properties:
\begin{enumerate}
\item let $\DD_{\sigma}(p_i)=\{D_\sigma(p_i,\alpha):\alpha\in\Lambda\}$, 
then for each $i$ and $\alpha$ there exists a sub-bundle $\DD(p_i,\alpha)\to\TT$ 
such that the fiber over $\sigma\in\TT$ is $D_\sigma(p_i,\alpha)$.
\item  each $\DD_{\sigma}(p_i)$ is invariant under the action of $\pi_1(S_{g,n})$,
\item elements of $\DD_\sigma(p_1)\cup...\cup\DD_\sigma(p_n)$ are disjoint with each other,
\end{enumerate}

For each path $\gamma:I\to\TT$, $1\leq i\leq n$ and $\alpha\in\Lambda$,
there exists a pull back bundle $\DD_\gamma(p_i,\alpha)\to I$ 
such that the fiber over $t\in I$ is $D_{\gamma(t)}(p_i,\alpha)$.
Let $\DD_\gamma$ denote the collection $\{\DD_\gamma(p_i,\alpha):1\leq i\leq n,\alpha\in\Lambda\}$.
Consider a  subset $\BB$ of the  modulli space $\MM=\TT/MCG(S_{g,n})$, a path $\gamma:I\to\TT$ is said to be $\BB$-\textit{cobounded},
 if the image of $\gamma$ under the projection $\TT\to\MM$ lies in $\BB$. We prove the following theorem:
\begin{thm}\label{main}
Let $I$ be a closed, connected interval of $\mathbb R$. For a compact subset $\BB$ of the  moduli space $\MM=\TT/MCG(S_{g,n})$
 and for every $\rho\geq 1,\delta\geq 0$
there exists $P\geq 0$ such that the following holds:\\
If $\gamma:I\to\TT$ is $\BB$-cobounded and $\rho$-Lipschitz path, and if $\HH_\gamma$ is strongly $\delta$-hyperbolic relative to the
collection $\DD_\gamma$, then there exists a geodesic $\xi:I\to\TT$ joining end points of $\gamma$ such that the Hausdorff distance
between $\gamma(I)$ and $\xi(I)$ is at most $P$.
\end{thm}
Note that the fibers $\HH_\sigma={\mathbb H}^2\times\sigma$ of $\HH\to\TT$ are (uniformly) strongly hyperbolic relative to the collections
$\DD_\sigma=\{D_\sigma(p_i,\alpha):1\leq i\leq n, \alpha\in\Lambda\}$ of horodisks. 
Hence the coned-off spaces $\EE(\HH_\sigma,\DD_\sigma)$, $\sigma\in\TT$, are (uniformly) hyperbolic metric spaces.
Thus for a path $\gamma:I\to\TT$, there exists a bundle $\PP\HH_\gamma\to I$ of coned-off hyperbolic metric spaces
with fiber  $\EE(\HH_{\gamma(t)},\DD_{\gamma(t)})$. $\PP\HH_\gamma$ is also obtained by partially
electrocuting  each element $\DD_\gamma(p_i,\alpha)$ of $\DD_\gamma$ to a hyperbolic space $\LL_\gamma(p_i,\alpha)$, where 
 $\LL_\gamma(p_i,\alpha)$ is the locus of  cone points obtained by coning $D_{\gamma(t)}(p_i,\alpha)$ for all $t\in I$.
By Lemma 2.8 of \cite{comb}, 
if $\HH_\gamma$ is strongly hyperbolic relative to the collection $\DD_\gamma$ then $\PP\HH_\gamma$ is a hyperbolic metric
space.  
\begin{defn}
Given $\kappa>1$, a natural number $n$, $A\geq 0$, a sequence of positive numbers $\{r_j:j\in J\}$,
where $J$ is a subinterval of set of integers $\mathbb Z$, is said to satisfy $(\kappa,n,A)$-flaring property
if $j-n,j+n\in J$ and if $r_j>A$ then $\max\{r_{j-n},r_{j+n}\}\geq\kappa r_j$. 
\end{defn}

A path $\alpha:J\to \PP\HH_\gamma$, where $J\subset I$, is  said to be $\lambda$-\textit{quasivertical} 
if it is  $\lambda$-Lipschitz and also a section. 
Let $d_{\hhat\sigma}$ denote the metric of the coned-off space $\EE(\HH_\sigma,\DD_\sigma)$. 
Since $\PP\HH_\gamma$ is a hyperbolic space, so we have the following flaring properties:\\
\begin{prop}\label{flare}(Theorem 4.7 of \cite{comb})
 With the notations as above, given $\lambda\geq 1$ there exist $\kappa>1$, an integer $n\geq 1$ and a number $A>0$
 such that the following holds: \\
Let $\alpha,\beta:J\to\PP\HH_\gamma$ be two $\lambda$-quasivertical paths, then the sequence
 $s_j=d_{\hhat{(\gamma(j))}}(\alpha(j),\beta(j))$, where $j\in J\cap\mathbb Z$, satisfies $(\kappa,n,A)$-flaring property.
\end{prop}
We refer to \cite{mcg} for the definitions of measured foliation $\MM\FF$ and measured geodesic lamination $\MM\GG\LL$ of general
hyperbolic surfaces. For each $\mu\in\MM\FF$, let $\mu_t$ denote the measured geodesic lamination on the hyperbolic
surface $\SSS_{\gamma(t)}=\HH_{\gamma(t)}/\pi_1(S_{g,n})$. Let $\SSS^b_{\gamma(t)}$ denote the `thick part' of $\SSS_{\gamma(t)}$
i.e. $\SSS^b_{\gamma(t)}$ is obtained from $\SSS_{\gamma(t)}$ by deleting the images of interior of horodisks under the projection 
$\HH_{\gamma(t)}\to\SSS_{\gamma(t)}$. Now each $\mu\in\MM\FF$ induce a geodesic lamination $\mu^b_t (\subset\mu_t)$ on $\SSS^b_{\gamma(t)}$.
A connection path of the sub-bundle $\SSS^b_{\gamma}\to I$ is a piecewise smooth section of the projection map
which is everywhere tangent to the connection on the bundle $\SSS^b_{\gamma}\to I$. The connection map 
$h_{st}:\SSS^b_{\gamma(s)}\to\SSS^b_{\gamma(t)}$ ($s\leq t$) is defined by moving points of $\SSS_{\gamma(s)}$
to $\SSS_{\gamma(t)}$ along connection paths. 
In \cite{farb}, it was proved that connection maps  $h_{st}$ are 
bilipschitz maps. 
For $\mu\in\MM\FF$ and $\sigma\in\TT$, the length of $\mu$ with respect
to $\sigma$ is defined by $len_\sigma(\mu)=\int d\mu$.
From proposition \ref{flare},  it follows that for any leaf segment $l_s$ of $\mu_s$, the sequence of
lengths $len_{s+i}(h_{s,s+i}(l_s))$ satisfies the flaring property.
As a consequence, we have the following theorem :
\begin{thm}\label{len flare}(Lemma 3.6 of \cite{mos})
For a compact subset $\BB$ of the  moduli space $\MM$
and for every $\rho\geq 1$, there exist constants $L\geq 1,\kappa>1, n\in\mathbb Z_+$ such that the following holds:
Let $\gamma:I\to\TT$ be a $\BB$-cobounded and $\rho$-Lipschitz path, for any $\mu\in\MM$,
the sequence $i\to len_{\gamma(i)}(\mu^b),~(i\in I\cap\mathbb Z)$, satisfies the $L$-Lipschitz,
$(\kappa,n,0)$-flaring property.
\end{thm}
For $\mu\in\MM\FF$, we say $\mu$ is \textit{realized} at  $p$, where $p$ is a finite number or $p\in\{-\infty,+\infty\}$,
if $len_{\gamma(i)}(\mu)$ achieves minimum at $p$.
\begin{prop}\label{realize}(Proposition 3.12 of \cite{mos})
For each $k\in I\cap\mathbb Z$, there exists $\mu\in\MM\FF$ which is finitely realized. If $I$ is infinite,
for each infinite end $\pm\infty$ of $I$ there exists $\mu_{\pm}\in\MM\FF$ which is realized at $\pm\infty$ respectively.
\end{prop}
Now for a compact subset $\BB\subset\MM$, numbers $\rho\geq 1,\delta\geq 0,\eta>0$, consider $\Gamma_{\beta,\rho,\delta,\eta}$
to be the set of all triples $(\gamma,\mu_-,\mu_+)$  with the following properties (see \cite{mos}):
\begin{enumerate}
\item $\gamma:I\to\TT$ is  $\BB$-cobounded, $\rho$-Lipschitz path,
such that $\HH_\gamma$ is $\delta$-hyperbolic relative to $\DD_\gamma$,
\item $0\in I$, and each $\mu_{\pm}\in\MM\FF$ is normalized to have length $1$ in the hyperbolic structure $\gamma(0)$,
\item the lamination $\mu_+$  is realized in $\SSS_{\gamma}$ near the right  end
in the following way:\\
(a) If $I$ is right  infinite, then $\mu_+$ is realized at $+\infty$,\\
(b) If $I$ is right finite, with right end point $M$, then there exists a minimum of
length sequence $len_{\gamma(i)}(\mu_+)$ lying in the interval $[M-\eta,M]$.\\
The lamination $\mu_-$ is realized similarly in $\SSS_{\gamma}$ near the left end.
\end{enumerate}
Let $\AAA\subset\TT$ be a compact set such that   each $(\gamma,\mu_-,\mu_+)\in \Gamma_{\beta,\rho,\delta,\eta}$,
may be translated by the action of $MCG(S_{g,n})$ so that $\gamma(0)\in\AAA$.
If $\gamma_i$ converges to $\gamma$, then  in the Gromov-Hausdorff topology,
  $\HH_{\gamma_i}$ converges to $\HH_{\gamma}$  and $\DD_{\gamma_i}$
converges to  $\DD_{\gamma}$. Hence, $\GG(\HH_{\gamma_i},\DD_{\gamma_i})$ converges to $\GG(\HH_\gamma,\DD_\gamma)$
in the Gromov-Hausdorff topology. The Gromov-Hausdorff limit of a sequence of $\delta$-hyperbolic spaces 
is $\delta$-hyperbolic (\cite{gro}). Therefore, if $\HH_{\gamma_i}$ are $\delta$-hyperbolic
relative to $\DD_{\gamma_i}$ for all $i$, then $\HH_{\gamma}$ is also $\delta$-hyperbolic relative to $\DD_{\gamma}$.
This justifies  the set
 $\AAA_{\beta,\rho,\delta,\eta}=\{(\gamma,\mu_-,\mu_+)\in\Gamma_{\beta,\rho,\delta,\eta}:\gamma(0)\in\AAA\}$ is compact.

\begin{prop}\label{cocpt}\cite{mos}
 The action of $MCG(S_{g,n})$ on $\Gamma_{\beta,\rho,\delta,\eta}$ is cocompact.
\end{prop}
\noindent{\textbf{Proof of Theorem \ref{main}}}\\
For $(\gamma,\mu_-,\mu_+)\in\Gamma_{\beta,\rho,\delta,\eta}$, 
let $a_-(t)=\frac{1}{len_{\gamma(t)}(\mu_-)}$ and $a_+(t)=\frac{1}{len_{\gamma(t)}(\mu_+)}$. $\mu_-,\mu_+$ fills $S_{g,n}$ (See \cite{mos}),
therefore $\mu_-,\mu_+$ defines a conformal structure $\sigma(\mu_-,\mu_+)$ on $S_{g,n}$.
Consider the map $\xi(t)=\sigma(a_-(t)\mu_-,a_+(t)\mu_+)$, $t\in I$, then the image of the map $\xi:I\to\TT$
is a geodesic in $\TT$ joining $\mu_-$ and $\mu_+$. For $i\in I\cap\mathbb Z$, define $\gamma'(s)=\gamma(s+i)$,
then the triple $(\gamma',a_-(t)\mu_-,a_+(t)\mu_+)$ lies in a translate of the compact set $\AAA_{\beta,\rho,\delta,\eta}$ 
by an element of $MCG(S_{g,n})$. The map taking $(\alpha,\lambda_-,\lambda_+)\in\Gamma_{\beta,\rho,\delta,\eta}$ to
$(\alpha(0),\sigma(\lambda_-,\lambda_+))\in\TT\times\TT$ is $MCG(S_{g,n})$ equivariant and continuous and hence
has $MCG(S_{g,n})$ cocompact image. Therefore, the Teichmuller distance $d_{\TT}$ between
 $\gamma(i)$ and $\sigma(a_-(i)\mu_-,a_+(i)\mu_+))=\xi(i)$ is bounded.
Now for $t\in I$, there exists $i\in I\cap\mathbb Z$ such that $|t-i|\leq 1$. 
As $\gamma$ is $\rho$-Lipschitz, therefore $d_{\TT}(\gamma(t),\gamma(i))\leq\rho$.
Also, there exists $L>0$ such that $d_{\TT}(\xi(t),\xi(i))\leq L$ (See \cite{mos}).
 Thus, there exists $P>0$ such that the Hausdorff distance between
$\gamma$ and $\xi$  is at most $P$. \qed
\section{Application}
Consider the following short exact sequence of pair of finitely generated groups:\\
$$1\to (\pi_1(S_{g,1}),K_1)\to (G,N_G(K_1))\to (Q,Q)\to 1,$$
where $K_1$ is  peripheral subgroup of $\pi_1(S_{g,1})$,
$G$ is strongly hyperbolic relative to $N_G(K_1)$ and $Q$ is a subgroup of $MCG(S_{g,1})$. 
Let $\Phi: Q\to \TT$ denote the orbit map, then for any geodesic $\gamma':I\to Q$, 
$\gamma=\Phi\circ\gamma ': I\to\TT$ is a cobounded and Lipschitz path. Since $G$ is strongly hyperbolic
relative to $N_G(K_1)$,  the bundle $\EE(G,K_1)$ over $Q$ is hyperbolic. Hence,
$\EE(G,K_1)\to Q$ satisfies flaring property. In particular, the sub-bundle $\PP\HH_{\gamma}\to I$
satisfies the flaring property. Therefore, by the converse of strong combination theorem in \cite{comb},
 $\HH_{\gamma}$ is strongly hyperbolic relative to $\DD_{\gamma}$. Hence, as an application of 
Theorem \ref{main}, $Q$ is a convex cocompact subgroup of $MCG(S_{g,1})$. The converse of this result
is also true (see \cite{mj}). So, we have the following theorem :
\begin{thm}\cite{mj}
Consider the following short exact sequence of pair of finitely generated groups
$$1\to (\pi_1(S_{g,1}),K_1)\to (G,N_G(K_1))\to (Q,Q)\to 1,$$
where $\pi_1(S_{g,1})$ is strongly hyperbolic
relative to $K_1$.  $G$ is strongly hyperbolic relative to $N_G(K_1)$ if and only if
$Q$ is a convex cocompact subgroup of $MCG(S_{g,1})$
\end{thm}

\end{document}